\documentclass[a4 paper,12pt,bezier]{article}
\usepackage[top=3cm,bottom=3cm,left=2.5cm,right=2.5cm]{geometry}
\usepackage{tikz, color,xcolor}
\usepackage{caption}
\usepackage{amsmath}
\usepackage{graphics}
\usepackage{mathrsfs}
\usepackage{subcaption}
\usepackage{amsmath}
\usepackage{amsfonts,amsthm,amssymb,mathrsfs,bbding}
\usepackage{float}
\usepackage{graphics,epstopdf}
\usepackage{graphics,multicol}
\usepackage{graphicx}
\usepackage{color}
\usepackage{enumerate}
\usepackage{caption}
\usepackage{rotating}
\usepackage{lscape}
\usepackage{longtable}
\allowdisplaybreaks[4]
\usepackage{tabularx}
\usepackage[colorlinks=true,anchorcolor=blue,filecolor=blue,linkcolor=blue,urlcolor=blue,citecolor=blue]{hyperref}
\usepackage{extarrows}
\usepackage{cite}
\usepackage{latexsym,bm}
\usepackage{mathtools}
\usepackage{changes}
\usepackage{lipsum}
\definechangesauthor[name={kk}, color=blue]{kk}
\newtheorem{theorem}{Theorem}
\newtheorem{prob}{Problem}

\newtheorem{lemma}{Lemma}

\newtheorem{claim}{Claim}

\renewcommand\proofname{\bf Proof}

\definecolor{ECUST}{RGB}{0,106,167}
\addtocounter{section}{0}

\makeatletter
\newsavebox\myboxA
\newsavebox\myboxB
\newlength\mylenA
\newcommand*\xoverline[2][0.75]{%
\sbox{\myboxA}{$\m@th#2$}%
\setbox\myboxB\null
\ht\myboxB=\ht\myboxA%
\dp\myboxB=\dp\myboxA%
\wd\myboxB=#1\wd\myboxA
\sbox\myboxB{$\m@th\overline{\copy\myboxB}$}
\setlength\mylenA{\the\wd\myboxA}
\addtolength\mylenA{-\the\wd\myboxB}%
\ifdim\wd\myboxB<\wd\myboxA%
\rlap{\hskip 0.5\mylenA\usebox\myboxB}{\usebox\myboxA}%
\else
\hskip -0.5\mylenA\rlap{\usebox\myboxA}{\hskip 0.5\mylenA\usebox\myboxB}%
\fi}
\makeatother
\begin{document}
\title{\bf The largest eigenvalue of $\mathcal{C}_4^{-}$-free signed graphs\footnote{This work is supported by the National Natural Science Foundation of China (Grant No. 12271162 ), Natural Science Foundation of Shanghai (No. 22ZR1416300) and The Program for Professor of Special Appointment (Eastern Scholar) at Shanghai Institutions of Higher Learning (No. TP2022031).}}
\author{Yongang Wang, Huiqiu Lin\footnote{Corresponding author. Email: huiqiulin@126.com}\\[2mm]
\small School of Mathematics, East China University of Science and Technology, \\
\small  Shanghai 200237, P.R. China\\}
\date{}
\maketitle
{\flushleft\large\bf Abstract}
 Let  $\mathcal{C}_{k}^{-}$ be the set of all   negative $C_k$. For odd cycle, Wang, Hou and Li \cite{C3free} gave 
 a spectral  condition for the existence of negative $C_3$ in  unbalanced signed graphs. 
 For even cycle,   we determine the maximum index among all $\mathcal{C}_4^{-}$-free unbalanced signed graphs  and completely characterize the extremal signed graph in this paper.  This could be regarded as a
signed graph  version of the results by Nikiforov\cite{NikiKr} and Zhai and Wang\cite{zhaiC4}.

\begin{flushleft}
\end{flushleft}

\section{Introduction}
All graphs in this paper are simple.  Let $\mathcal{F}$ be a family of graphs. A graph $G$ is $\mathcal{F}$-free if $G$ does not contain any graph in $\mathcal{F}$
as a subgraph. The classical spectral Tur\'{a}n  problem is to determine the maximum spectral radius of an  $\mathcal{F}$-free graph of order $n$, which is known as  the \textit{spectral Tur\'{a}n number} of $\mathcal{F}$. 
This problem was originally proposed by Nikiforov\cite{proposeNikiforov}. 
With regard to  unsigned graphs, much attention has been paid to the spectral Tur\'{a}n  problem in the past decades, see \cite{1Babai,triangleLin,WilfKr,NikiC4,Li,JCTBlin,Bollobas}. In this paper, we focus on the spectral Tur\'{a}n  problem in signed graphs.

 A  \textit{signed graph}   $\Gamma=(G,\sigma)$  consists of a graph $G=(V(G),E(G))$ and a sign function $\sigma : E\rightarrow \{-1,1\}$, where $G$ is its underlying graph and  $\sigma$ is its sign function. An edge $e$ is \textit{positive}\ (\textit{negative})\ if $\sigma(e)=1$ (resp.  $\sigma(e)=-1$).  A cycle  $C$ in  a signed graph $\Gamma$ is called \emph{positive} (resp.\ \emph{negative})  if the number of its negative edges is even (resp.\ odd). A signed graph is called \emph{balanced} if  all  its cycles are positive; otherwise, it is called \emph{unbalanced}.  The \textit{adjacency} \textit{matrix} of  $\Gamma$ is denoted by $A(\Gamma)=(a^{\sigma}_{ij})$,\ where $a^{\sigma}_{ij} =\sigma(v_{i}v_{j})$ if $v_{i}\sim v_{j}$,\ and $0$ otherwise. The eigenvalues of $A({\Gamma})$ are called the eigenvalues of  $\Gamma$.   The largest eigenvalue of ${\Gamma}$ is called the \textit{index} of $\Gamma$ and denoted by $\lambda_1(\Gamma)$. For more details about the notion of signed graphs, we refer to\cite{treepositive,cycle}.

The spectral Tur\'{a}n problem of signed graphs  has been studied in recent years.  Let $\mathcal{K}_4^{-}$ be the set of all  unbalanced $K_4$.   Chen and Yuan\cite{K4free} gave the spectral Tur\'{a}n number of $\mathcal{K}_4^{-}.$   For the largest eigenvalue  of a signed graph with   certain structures, 
Koledin and Stani\'{c}\cite{connectedindex} studied connected signed graphs of fixed order, size and number
of negative edges that maximize the index  of their adjacency matrices. After that, signed graphs maximizing
the index in suitable subsets of  signed complete graphs have been studied by Ghorbani and Majidi\cite{maxindex},  
Li, Lin and Meng\cite{Meng} and   Akbari,  Dalvandi,  Heydari and  Maghasedi\cite{kedge}.  
It is well known that the eigenvalues of a balanced signed graph are the same as those of its underlying graph. Therefore, the largest eigenvalue  of an unbalanced signed graph  has attracted more attention of scholars.   In  $2019$, Akbari, Belardo, Heydari, Maghasedi and Souri\cite{unicyclic}  determined the 
 signed graphs achieving the minimal or  maximal index in the class of unbalanced signed unicyclic graphs. In 2021,   He, Li, Shan and Wang \cite{shuangquan} gave the first five largest indices among
 all unbalanced signed bicyclic graphs of order  $n\geq36$.    In $2022$,  Brunetti and Stani\'{c}\cite{unbalancedindex} studied the extremal spectral radius among  all unbalanced connected signed graphs.   More results on the  the spectral theory of signed graphs can be found in \cite{open,signedturan,WYQ,JCTA,huanghao}, where \cite{open} is an excellent survey about some general results and  problems on the spectra of signed graphs.  

 The study of cycles from the eigenvalue perspective has a long history, such as $C_{2k+1}$ \cite{oddNiki}, $C_4$ \cite{NikiKr,zhaiC4}, $C_6$ \cite{C6Zhai},  $C_{2k}$ for $k\geq 4$ \cite{CDT}, cycles of consecutive lengths \cite{LiNing,NP,ZL,oddNiki} and long cycles \cite{LiEJC,GH}.  Let  $\mathcal{C}_{k}^{-}$ be the set of all   negative $C_k$.   For signed graphs, Wang, Hou and Li \cite{C3free} determined the spectral Tur\'{a}n number of $\mathcal{C}_3^{-}.$   Denote by $(G,+)$ (resp. $(G,-)$) the signed graph whose edges are all positive (resp. negative).
Note that $(K_n,-)$ is unbalanced and $\mathcal{C}_{2k}^{-}$-free  whose spectral radius  is always $n-1$ for $n\geq4$.
Then it is interesting to study the existence  of negative $C_{2k}$ from the largest eigenvalue condition (see Problem \ref{problem}).
\begin{prob}\label{problem}
What is the largest eigenvalue among all $\mathcal{C}_{2k}^{-}$-free  unbalanced  signed graphs for $k\geq 2?$
\end{prob}

\begin{figure}[htbp!]
  \centering
\begin{subfigure}[c]{0.36\linewidth}
    \centering
    \begin{tikzpicture}
      \tikzstyle{every node}=[font=\small]

\node (4) at (0,0) [shape=circle,line width=1.26pt,draw,fill=gray!25,inner sep=10pt,minimum size=19.7mm] {$K_{n-2}$};
       \node (3) at (-1,0) [shape=circle,draw,fill=black!100,inner sep=2.1pt,minimum size=2mm] {};
      \node (1) at (-2.4,-1.1) [shape=circle,draw,fill=black!100,inner sep=1pt,minimum size=2mm] {};
      \node (2) at (-2.4,1.1) [shape=circle,draw,fill=black!100,inner sep=1pt,minimum size=2mm] {};
      \node (5) at (-1.16,0.45) [] {$v_3$};
      \node (6) at (-2,1.1) [] {$v_1$};
      \node (7) at (-2,-1.1) [] {$v_2$};
      \draw [red,line width=1.2pt] (1) to (2);
       \draw [black,line width=1.2pt] (3) to (1);
      \draw [black,line width=1.2pt] (3) to (2);
  \end{tikzpicture}
\captionsetup{labelformat=empty}
    \subcaption{$\Gamma_1$}
\end{subfigure}
 \begin{subfigure}[c]{0.36\linewidth}
    \centering
    \begin{tikzpicture}
      \tikzstyle{every node}=[font=\small]
      \node (0) at (0,0) [shape=circle,draw,fill=black!100,inner sep=1pt,minimum size=2mm] {};
      \node (1) at (-0.75,1.15) [shape=circle,draw,fill=black!100,inner sep=1pt,minimum size=2mm] {};
      \node (2) at (-0.75,-1.15) [shape=circle,draw,fill=black!100,inner sep=1pt,minimum size=2mm] {};
      \node (3) at (-1.5,0) [shape=circle,draw,fill=black!100,inner sep=1pt,minimum size=2mm] {};
      \node (4) at (2,0) [shape=circle,draw,fill=gray!25,line width=1.26pt,inner sep=1pt,minimum size=19.3mm] {$K_{n-4}$};
      \node (5) at (0.4,0)[]{$v_2$};
       \node (6) at (-1.9,0)[]{$v_1$};
       \node (7) at (-1.2,1.15)[]{$v_3$};
       \node (8) at (-1.2,-1.15)[]{$v_4$};
      \draw [black,line width=1.2pt] (0) to (1);
      \draw [red,line width=1.2pt] (0) to (3);
       \draw [black,line width=1.2pt] (0) to (2);
      \draw [green,line width=1.2pt] (2) to (4);
      \draw [green,line width=1.2pt] (1) to (4);
      \draw [black,line width=1.2pt] (2) to (3);
      \draw [black,line width=1.2pt] (3) to (1);
  \end{tikzpicture}
\captionsetup{labelformat=empty}
    \subcaption{$\Gamma_2$}
\end{subfigure}
	\caption{The signed graphs $\Gamma_1$ and $\Gamma_2$.}\label{fig2}
\end{figure}
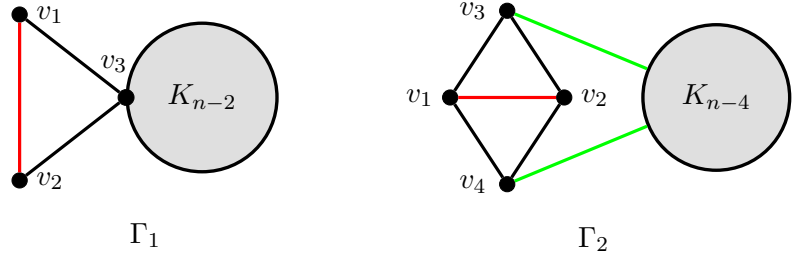

Suppose $\Gamma=(G,\sigma)$ is a signed graph and $U\subset V(G)$.
The operation that changes the sign of all edges between $U$ and $V(G)\backslash U$ is called a \emph{switching operation}. If a signed graph $\Gamma^{\prime}$ is obtained from $\Gamma$ by applying finitely many switching operations, then $\Gamma$ is said to be \emph{switching equivalent} to $\Gamma^{\prime}$.  For $n\geq5$, we define the signed graphs $\Gamma_1$ and $\Gamma_2$ as shown in Figure \ref{fig2}, where the circles represent unsigned complete graphs, the red lines represent negative edges and the other lines represent positive edges, and especially,  the green lines joined to the circles represent the connection of all possible edges. In this paper, we give an answer to Problem \ref{problem} for $k=2$ as follows.

\begin{theorem}\label{C4-free}
Let $\Gamma=(G,\sigma)$ be an unbalanced signed graph of order $n\geq5$. If $$\lambda_{1}(\Gamma)\geq\lambda_{1}(\Gamma_{1}),$$  then  $\Gamma$ contains a negative $C_4$ unless $\Gamma$ is switching equivalent to $\Gamma_{1}$ (see Figure \ref{fig2}). 
\end{theorem}

\section{The largest eigenvalues of signed graphs  $\Gamma_1$ and $\Gamma_2$}
In this section, we shall show that $\lambda_{1}(\Gamma_1)>n-3.$  We now introduce the definition of equitable  quotient matrix.

Let $M$ be a real symmetric matrix of order $n$, and let $[n]=\{1,2, \ldots, n\}$. Given a partition $\Pi:[n]=X_1 \cup X_2 \cup \cdots \cup X_k$, the matrix $M$ can be written as
$$
M=\left[\begin{array}{cccc}
M_{1,1} & M_{1,2} & \cdots & M_{1, k} \\
M_{2,1} & M_{2,2} & \cdots & M_{2, k} \\
\vdots & \vdots & \ddots & \vdots \\
M_{k, 1} & M_{k, 2} & \cdots & M_{k, k}
\end{array}\right] .
$$
If all row sums of $M_{i, j}$ are the same, say $b_{i, j}$, for all $i, j \in\{1,2, \ldots, k\}$, then $\Pi$ is called an $equitable$   $partition$ of $M$, and the matrix $Q=\left(b_{i, j}\right)_{i, j=1}^k$ is called an $equitable$  $quotient$  $matrix$ of $M$.

\begin{lemma}{\bf(\cite[p.24]{fz})}\label{quotient}
Let $M$ be a real symmetric matrix, and let $Q$ be an equitable quotient matrix of $M$. Then the eigenvalues of $Q$ are also eigenvalues of $M$.
\end{lemma}

\begin{lemma}\label{dayu n3}
Let $\Gamma_1$ and $\Gamma_2$ be the signed graphs as shown in Figure \ref{fig2}. Then we have the following statements.
\begin{itemize}
\item[(1)] $\lambda_{1}(\Gamma_1)>n-3.$
\item[(2)] $\lambda_{1}(\Gamma_2)< \lambda_{1}(\Gamma_1).$
\end{itemize}
\end{lemma}
\renewcommand\proofname{\bf Proof}
\begin{proof}
$(1)$ Let  $J$, $I$  and $O$ denote the  all-ones matrix, identity matrix and all-zeros matrix, respectively. By a suitable partition,
$$A(\Gamma_1)=\left[\begin{array}{cccc}
     0 & -1 & 1 & O\\
     -1 & 0& 1  & O\\
     1 & 1& 0 & J\\
     O & O& J& J-I\\
\end{array}\right],$$
and $A(\Gamma_1)$ has the equitable quotient matrix 
$$Q_1=\left[\begin{array}{cccc}
     0 & -1 & 1 & 0\\
     -1 & 0& 1  & 0\\
     1 & 1& 0 & n-3\\
     0 & 0& 1& n-4\\
\end{array}\right].$$
Note that rank$(A(\Gamma_1)+I)$=4. Then  $A(\Gamma_1)+I$ has $0$ as an eigenvalue with multiplicity $n-4.$ Hence, $A(\Gamma_1)$ has $-1$ as an eigenvalue with multiplicity $n-4.$
By a simple calculation, the characteristic polynomial of $Q_1$ is
$$f(x)=(x-1)(x^3+(5-n)x^2+(5-2n)x+n-5).$$
It is easy to check that $f(n-3)=-2n+8<0$ and $f(n-2)=(n-3)^2(n+1)>0$. Then $\lambda_1(Q_1)>n-3.$ Since $f(-1)\neq0,$ we have  $\lambda_1(\Gamma_1)=\lambda_1(Q_1)>n-3$ by Lemma \ref{quotient}.


$(2)$ For $5\leq n\leq 6,$ by a direct calculation, we have $\lambda_{1}(\Gamma_2)< \lambda_{1}(\Gamma_1).$ For $n\geq7,$  by a suitable partition,
$$A(\Gamma_2)=\left[\begin{array}{cccc}
     0 & -1 & J & O\\
     -1 & 0& J  & O\\
     J & J& O_{2\times 2} & J\\
     O & O& J& J-I\\
\end{array}\right],$$
and $A(\Gamma_2)$ has the equitable quotient matrix 
$$Q_2=\left[\begin{array}{cccc}
     0 & -1 & 2 & 0\\
     -1 & 0& 2  & 0\\
     1 & 1& 0 & n-4\\
     0 & 0& 2& n-5\\
\end{array}\right].$$
Note that rank$(A(\Gamma_2))\leq n-1$ and  rank$(A(\Gamma_2)+I)=5$. Then $A(\Gamma_2)$ has $0$ as an eigenvalue with multiplicity at least  $1,$ and  $A(\Gamma_2)$ has $-1$ as an eigenvalue with multiplicity $n-5.$
By a simple calculation, the characteristic polynomial of $Q_2$ is
$$g(x)=(x-1)h(x),$$  where  $h(x)=x^3+(6-n)x^2+(9-3n)x+2n-12$.
Observe that  $h(-\infty)<0,$ $h(0)>0,$  $h(n-4)<0$ and  $h(n-3)>0.$ Then the three roots of $h(x)$ lie in $(-\infty,0),$  $(0,n-4)$ and $(n-4,n-3)$, respectively. Since $g(0)\neq0$ and $g(-1)\neq0,$  we have  $\lambda_1(\Gamma_2)=\lambda_1(Q_2)<n-3<\lambda_1(\Gamma_1)$ by Lemma \ref{quotient} and Lemma \ref{dayu n3}$(1)$.
\end{proof}

\section{Proof of Theorem \ref{C4-free}}
By the table of the spectra of signed graphs with five vertices\cite{table}, we can check that Theorem \ref{C4-free} is true for $n=5.$ Therefore, we now assume that $n\geq 6.$  We first give two lemmas  which are needed  in the proof of Theorem \ref{C4-free}.

\begin{lemma}{\bf  ( \cite[Lemma 2.5]{SGX})}\label{Nonnegative}
Let $\Gamma$ be a signed graph. Then there exists a signed graph $\Gamma^{\prime}$ switching equivalent to $\Gamma$ such that   $A(\Gamma^{\prime})$ has a non-negative eigenvector corresponding to $\lambda_{1}(\Gamma^{\prime})$.
\end{lemma}

\begin{lemma}{\bf (\cite[Proposition 3.2]{treepositive})}\label{remaincycle}
 Two signed graphs with the same underlying graph are switching equivalent
 if and only if they have the same set of positive cycles.
\end{lemma}

Nikiforov\cite{NikiKr} and Zhai and Wang\cite{zhaiC4} determined the spectral conditions for the existence of $C_4$  for odd $n$ and even $n$, respectively, and they also further characterized the corresponding    spectral extremal graphs.  
 \begin{lemma}\label{GC4free}
 Let $G$ be a $C_4$-free graph of order $n$ with  $\lambda_1(G)=\lambda.$  Then  
 
     $(i)${\bf (\cite{NikiKr})} If $n$ is odd, then $\lambda^{2}-\lambda-(n-1) \leq  0.$ 
     
     $(ii)${\bf (\cite{zhaiC4})} If $n$ is even, then  $\lambda^{3}-\lambda^{2}-(n-1)\lambda+1 \leq  0.$
 \end{lemma}

Now, we  are in a position to give the proof of Theorem \ref{C4-free}.
\renewcommand\proofname{\bf Proof of Theorem \ref{C4-free}}
\begin{proof}
Suppose that $\Gamma=(G,\sigma)$ has the maximum index among all   $\mathcal{C}_{4}^{-}$-free unbalanced signed graphs. We  shall show that $\Gamma$ is switching equivalent to $\Gamma_1.$  Let  $\Gamma^{\prime}=(G,\sigma^{\prime})$ be a signed graph switching equivalent to $\Gamma.$ By Lemma \ref{Nonnegative}, we can assume that  $A(\Gamma^{\prime})$ has a non-negative eigenvector corresponding to $\lambda_{1}(\Gamma^{\prime})$.   Then by Lemma \ref{remaincycle},     $\Gamma^{\prime}$ is also unbalanced and $\mathcal{C}_{4}^{-}$-free. Furthermore,  $\Gamma^{\prime}$ also has the maximum index among all $\mathcal{C}_{4}^{-}$-free unbalanced  signed graphs. Set $V(\Gamma^{\prime})=\{v_1,v_2,\ldots,v_n\}.$
Let $x=(x_1,x_2,\ldots,x_n)^{\top}$ be the non-negative unit eigenvector of $A(\Gamma^{\prime})$ corresponding to  $\lambda_1(\Gamma^{\prime}),$ where $x_i$ corresponds to the vertex $v_i$ for $1\leq i\leq n.$  Then $$\lambda_1(\Gamma^{\prime})=x^{\top}A(\Gamma^{\prime})x.$$
Since $\Gamma_{1}$ is unbalanced and $\mathcal{C}_{4}^{-}$-free, we may suppose that  $\lambda_1(\Gamma^{\prime})\geq\lambda_1(\Gamma_{1})>n-3$ by Lemma \ref{dayu n3}. Now we begin to analyze the structure of $\Gamma^{\prime}$. First we give  some claims.
\begin{claim}\label{1zero}\label{Claim1}
$x$ has at most one zero coordinate.
\end{claim}
Otherwise, without loss of generality, assume $x_1=x_2=0,$ then $$\begin{aligned}
\lambda_1(\Gamma^{\prime})&=x^{\top}A(\Gamma^{\prime})x=(x_3,x_4,\ldots,x_n)A(\Gamma^{\prime}-v_1-v_2)(x_3,x_4,\ldots,x_n)^{\top}\\
&\leq\lambda_1(\Gamma^{\prime}-v_1-v_2)\leq\lambda_1(K_{n-2})=n-3,
\end{aligned} $$
a contradiction. So Claim \ref{Claim1} holds.

Let $N_{\Gamma^{\prime}}(v_i)$ denote the  set of neighbours of  $v_i$ in $\Gamma^{\prime}.$ 
\begin{claim}\label{connected}
  $\Gamma^{\prime}$ is connected.
\end{claim}
Otherwise, assume $\Gamma^{\prime}_1$ and $\Gamma^{\prime}_2$ are two distinct connected components of $\Gamma^{\prime}$, where $\lambda_1(\Gamma^{\prime})=\lambda_1(\Gamma^{\prime}_1)$.  Without loss of generality, we choose two vertices $v_i\in V(\Gamma^{\prime}_1)$ and $v_j\in V({\Gamma^{\prime}_2}).$  Then we can construct a new signed graph $\Gamma^{\ast}$   obtained from $\Gamma^{\prime}$ by adding a positive edge $v_iv_j.$  Clearly, $\Gamma^{\ast}$ is also unbalanced and $\mathcal{C}_4^-$-free. By Rayleigh principle, we obtain that
$$\begin{aligned}
\lambda_1(\Gamma^{\ast})-\lambda_1(\Gamma^{\prime})&\geq x^{\top}A(\Gamma^{\ast})x-x^{\top}A(\Gamma^{\prime})x\\
&=2x_ix_j\geq0.
\end{aligned} $$
If $\lambda_1(\Gamma^{\ast})=\lambda_1(\Gamma^{\prime})$, then $x$ is also an eigenvector of $A(\Gamma^{\ast})$
corresponding to $\lambda_1(\Gamma^{\ast}).$ Based on the following equations,
$$\lambda_1(\Gamma^{\prime})x_{i}=\sum_{v_s\in N_{\Gamma^{\prime}}(v_i)}\sigma^{\prime}(v_sv_i)x_s,$$
$$\lambda_1(\Gamma^{\prime})x_{j}=\sum_{v_s\in N_{\Gamma^{\prime}}(v_j)}\sigma^{\prime}(v_sv_j)x_s,$$
$$  \lambda_1(\Gamma^{\ast})x_{i}=\sum_{v_s\in N_{\Gamma^{\prime}}(v_i)}\sigma^{\prime}(v_sv_i)x_s+x_j$$ and
$$  \lambda_1(\Gamma^{\ast})x_{j}=\sum_{v_s\in N_{\Gamma^{\prime}}(v_j)}\sigma^{\prime}(v_sv_j)x_s+x_i,$$
we obtain that $x_i=x_j=0,$ which  contradicts Claim \ref{1zero}. Hence, $\lambda_1(\Gamma^{\ast})>\lambda_1(\Gamma^{\prime}),$ a contradiction.

Since $\Gamma^{\prime}$ is unbalanced,   $\Gamma^{\prime}$ contains at least one  negative edge and at least one negative cycle. Let $\mathscr{C}$ be one of the shortest negative cycles of $\Gamma^{\prime}$. 

\begin{claim}\label{containallnegative}
$\mathscr{C}$ contains all negative edges of  $\Gamma^{\prime}.$
\end{claim}
Otherwise,  without loss of generality,  assume $e=v_iv_j$ is a negative edge of $\Gamma^{\prime}$ and $e\notin E(\mathscr{C}).$ Then we can construct a new signed graph $\Gamma^{\ast}$  obtained from $\Gamma^{\prime}$ by deleting  $e.$  Clearly, $\Gamma^{\ast}$ is also unbalanced and $\mathcal{C}_4^-$-free. By Rayleigh principle, we obtain that
$$\begin{aligned}
\lambda_1(\Gamma^{\ast})-\lambda_1(\Gamma^{\prime})&\geq x^{\top}A(\Gamma^{\ast})x-x^{\top}A(\Gamma^{\prime})x\\
&=2x_ix_j\geq0
\end{aligned} $$
If $\lambda_1(\Gamma^{\ast})=\lambda_1(\Gamma^{\prime})$, then $x$ is also an eigenvector of $A(\Gamma^{\ast})$
corresponding to $\lambda_1(\Gamma^{\ast}).$ Based on the following equations,
$$\lambda_1(\Gamma^{\prime})x_{i}=\sum_{v_s\in N_{\Gamma^{\prime}}(v_i)}\sigma^{\prime}(v_sv_i)x_s,$$
$$\lambda_1(\Gamma^{\prime})x_{j}=\sum_{v_s\in N_{\Gamma^{\prime}}(v_j)}\sigma^{\prime}(v_sv_j)x_s,$$
$$  \lambda_1(\Gamma^{\ast})x_{i}=\sum_{v_s\in N_{\Gamma^{\prime}}(v_i)}\sigma^{\prime}(v_sv_i)x_s+x_j$$ and
$$  \lambda_1(\Gamma^{\ast})x_{j}=\sum_{v_s\in N_{\Gamma^{\prime}}(v_j)}\sigma^{\prime}(v_sv_j)x_s+x_i,$$
we obtain that $x_i=x_j=0,$  which  contradicts  Claim \ref{1zero}. Hence, $\lambda_1(\Gamma^{\ast})>\lambda_1(\Gamma^{\prime}),$ a contradiction.

\begin{claim}\label{youC4}
$G$ contains  $C_4$ as a subgraph.
\end{claim}
Otherwise, assume $G$ is $C_4$-free, then by Lemma \ref{GC4free},  $\lambda_1(\Gamma^{\prime})\leq \lambda_1(G)<n-3,$ a contradiction.

 Let $l=|V(\mathscr{C})|.$ Without loss of generality,  we can suppose   $\mathscr{C}=v_1v_2\cdots v_{l-1}v_lv_1.$ For $1\leq i<j\leq l,$ we assert that $v_{i}\nsim v_j$ if $j-i\geq2$ and $v_iv_j\neq v_1v_l$. Otherwise, there exists a shorter negative cycle than $\mathscr{C}$,  which  contradicts   the choice of $\mathscr{C}.$  Now, we define the  signed graphs  $\Gamma_{3}$ and $\Gamma_{4}$  as shown in Figure \ref{fig45}, where the black and red lines represent positive and negative edges,  respectively, and  especially,  the blue lines represent  the edges with uncertain signs.

\begin{figure}[htbp!]
  \centering
\begin{subfigure}[c]{0.36\linewidth}
    \centering
    \begin{tikzpicture}
      \tikzstyle{every node}=[font=\small]
       \node (0) at (0,0.3) [shape=circle,draw,fill=black!100,inner sep=2.1pt,minimum size=2mm] {};
      \node (1) at (1.5,0.3) [shape=circle,draw,fill=black!100,inner sep=1pt,minimum size=2mm] {};
      \node (2) at (-1,1.5) [shape=circle,draw,fill=black!100,inner sep=1pt,minimum size=2mm] {};
      \node (3) at (2.5,1.5) [shape=circle,draw,fill=black!100,inner sep=1pt,minimum size=2mm] {};
      \node (4) at (-0.8,-1.06) [shape=circle,draw,fill=black!100,inner sep=1pt,minimum size=2mm] {};
      \node (5) at (0.22,1.5) [shape=circle,draw,fill=black!100,inner sep=1pt,minimum size=2mm] {};
      \node (6) at (1.32,1.5) [shape=circle,draw,fill=black!100,inner sep=1pt,minimum size=2mm] {};
      \node (7) at (0.5,1.5) [shape=circle,draw,fill=black!100,inner sep=1pt,minimum size=0.1mm] {};
      \node (8) at (0.78,1.5) [shape=circle,draw,fill=black!100,inner sep=1pt,minimum size=0.1mm] {};
      \node (9) at (1.06,1.5) [shape=circle,draw,fill=black!100,inner sep=1pt,minimum size=0.1mm] {};
      \node (10) at (0.8,0.9) [] {$\mathscr{C}$};
      \node (11) at (-1.4,1.5) [] {$v_1$};
      \node (12) at (0.2,0.5) [] {$v_2$};
      \node (13) at (1.92,0.3) [] {$v_3$};
      \node (14) at (-0.3,-0.3) [] {$C_4^{\prime}$};
      \draw [red,line width=1.2pt] (0) to (2);
      \draw [red,line width=1.2pt] (0) to (1);
      \draw [blue,line width=1.2pt] (1) to (3);
       \draw [blue,line width=1.2pt] (2) to (5);
      \draw [blue,line width=1.2pt] (6) to (3);
      \draw [black,line width=1.2pt] (2) to (4);
      \draw [black,line width=1.2pt] (1) to (4);
  \end{tikzpicture}
\captionsetup{labelformat=empty}
    \subcaption{$\Gamma_3(l\geq5)$}
\end{subfigure}
  \begin{subfigure}[c]{0.26\linewidth}
    \centering
    \begin{tikzpicture}
      \tikzstyle{every node}=[font=\small]
       \node (4) at (1.5,2.36) [shape=circle,draw,fill=black!100,inner sep=2.1pt,minimum size=2mm] {};
      \node (1) at (0,-0.3) [shape=circle,draw,fill=black!100,inner sep=1pt,minimum size=2mm] {};
      \node (2) at (3,-0.3) [shape=circle,draw,fill=black!100,inner sep=1pt,minimum size=2mm] {};
      \node (3) at (1.5,1.04) [shape=circle,draw,fill=black!100,inner sep=1pt,minimum size=2mm] {};
      \node (5) at (1.5,0.12) [] {$\mathscr{C}$};
      \node (6) at (-0.2,0.01) [] {$v_1$};
      \node (7) at (1.5,0.7) [] {$v_2$};
      \node (8) at (3.15,0.06) [] {$v_3$};
       \node (9) at (1.5,1.5) [] {$C_4^{\prime}$};
      \draw [red,line width=1.2pt] (3) to (1);
      \draw [red,line width=1.2pt] (3) to (2);
      \draw [red,line width=1.2pt] (2) to (1);
      \draw [black,line width=1.2pt] (2) to (4);
      \draw [black,line width=1.2pt] (1) to (4);
  \end{tikzpicture}
\captionsetup{labelformat=empty}
    \subcaption{$\Gamma_4(l=3)$}
\end{subfigure}
	\caption{The signed graphs $\Gamma_3$ and $\Gamma_4$ in the proof of  Claim \ref{ALLedgeC4positive}.}\label{fig45}
\end{figure}
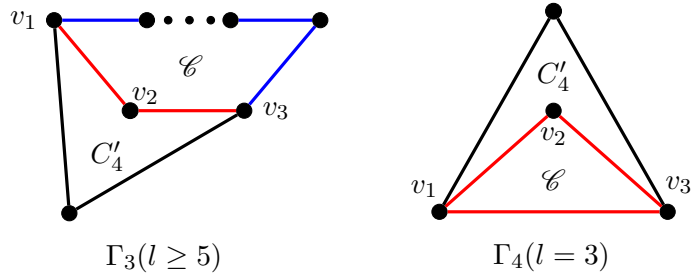

\begin{claim}\label{ALLedgeC4positive}
All edges of any cycle of $order$ $4$ in $\Gamma^{\prime}$ are positive.
\end{claim}
Otherwise,   assume that $C_4^{\prime}$ is a signed cycle of order $4$ in $\Gamma^{\prime}$ and   $C_4^{\prime}$ contains at least one  negative edge. Since  $\Gamma^{\prime}$ is $\mathcal{C}_4^{-}$-free,  $C_4^{\prime}$  contains two or four negative edges.    If $C_4^{\prime}$  contains four negative edges,  then by Claim \ref{containallnegative},   $\mathscr{C}$ contains four negative edges, and thus $\mathscr{C}$ is  positive,  a  contradiction.    If $C_4^{\prime}$  contains two negative edges, say $e_1$ and $e_2$,  then we assert that  $e_1$ and  $e_2$   must contain a common vertex. Otherwise, without loss of generality, let $e_1=v_1v_2$ and $e_2=v_3v_4.$    Then $e_1,e_2\in E(\mathscr{C})$ and $l\geq 5.$ 
Thus,  there  exists a shorter negative cycle than $\mathscr{C}$,  a contradiction. 
Without loss of generality, let $e_1=v_1v_2$ and $e_2=v_2v_3.$ Then $\Gamma^{\prime}$ must contain one of the signed graphs $\Gamma_{3}$ and  $\Gamma_{4}$ (see Figure \ref{fig45}) as a  subgraph, and we can construct a new signed graph $\Gamma^{\ast}$  obtained from $\Gamma^{\prime}$ by deleting $e_1$ and   $e_2.$  Clearly, $\Gamma^{\ast}$ is also unbalanced and $\mathcal{C}_4^-$-free. By Rayleigh principle, we obtain that
$$\begin{aligned}
\lambda_1(\Gamma^{\ast})-\lambda_1(\Gamma^{\prime})&\geq x^{\top}A(\Gamma^{\ast})x-x^{\top}A(\Gamma^{\prime})x\\
&=2x_2(x_1+x_3)\geq0.
\end{aligned} $$
If $\lambda_1(\Gamma^{\ast})=\lambda_1(\Gamma^{\prime})$, then $x$ is also an eigenvector of $A(\Gamma^{\ast})$
corresponding to $\lambda_1(\Gamma^{\ast}).$ Based on the following equations,
$$\lambda_1(\Gamma^{\prime})x_{2}=\sum_{v_s\in N_{\Gamma^{\prime}}(v_2)}\sigma^{\prime}(v_sv_2)x_s$$
and
$$  \lambda_1(\Gamma^{\ast})x_{2}=\sum_{v_s\in N_{\Gamma^{\prime}}(v_2)}\sigma^{\prime}(v_sv_2)x_s+x_1+x_3,$$
we obtain that $x_1=x_3=0,$ which contradicts  Claim \ref{1zero}. Hence, $\lambda_1(\Gamma^{\ast})>\lambda_1(\Gamma^{\prime}),$ a contradiction.
\begin{claim}\label{only1negative}\
$\Gamma^{\prime}$ contains exactly one negative edge.
\end{claim}
Otherwise, assume $\Gamma^{\prime}$ contains $m$ $(m\neq1)$ negative edges. By Claim \ref{containallnegative}, $m\geq3$ and $m$ is odd.  We first consider that $l=3.$ Without loss of generality, let $v_1v_2$ and $v_2v_3$ be two negative edges of $\Gamma^{\prime}.$ Then  we can construct a new signed graph $\Gamma^{\ast}$  obtained from $\Gamma^{\prime}$ by reversing the sign of   $v_1v_2$ and $v_2v_3.$  By Claim \ref{ALLedgeC4positive}, $\Gamma^{\ast}$ is also unbalanced and $\mathcal{C}_4^-$-free.  Furthermore, by Rayleigh principle, we obtain that
$$\begin{aligned}
\lambda_1(\Gamma^{\ast})-\lambda_1(\Gamma^{\prime})&\geq x^{\top}A(\Gamma^{\ast})x-x^{\top}A(\Gamma^{\prime})x\\
&=2(x_1x_2+x_2x_3)-2(-x_1x_2-x_2x_3)\\
&=4x_2(x_1+x_3)\geq0.
\end{aligned} $$
If $\lambda_1(\Gamma^{\ast})=\lambda_1(\Gamma^{\prime})$, then $x$ is also an eigenvector of $A(\Gamma^{\ast})$
corresponding to $\lambda_1(\Gamma^{\ast}).$
Based on the following equations,
$$\lambda_1(\Gamma^{\prime})x_{2}=\sum_{v_s\in N_{\Gamma^{\prime}}(v_2)}\sigma^{\prime}(v_sv_2)x_s$$ and
$$  \lambda_1(\Gamma^{\ast})x_{2}=\sum_{v_s\in N_{\Gamma^{\prime}}(v_2)}\sigma^{\prime}(v_sv_2)x_s+2(x_1+x_3),$$
we obtain that $x_1=x_3=0,$ which contradicts   Claim \ref{1zero}. Hence, $\lambda_1(\Gamma^{\ast})>\lambda_1(\Gamma^{\prime}),$ a contradiction. Now, we consider that  $l\geq 5.$  Without loss of generality, let  $v_1v_2$ and $v_3v_4$ be two negative edges of $\Gamma^{\prime}.$ Then  we can construct a new signed graph $\Gamma^{\ast}$  obtained from $\Gamma^{\prime}$ by reversing the sign of  $v_1v_2$ and $v_3v_4.$  By Claim \ref{ALLedgeC4positive}, $\Gamma^{\ast}$ is also unbalanced and $\mathcal{C}_4^-$-free. Furthermore, by Rayleigh principle and Claim \ref{1zero}, we obtain that
$$\begin{aligned}
\lambda_1(\Gamma^{\ast})-\lambda_1(\Gamma^{\prime})&\geq x^{\top}A(\Gamma^{\ast})x-x^{\top}A(\Gamma^{\prime})x\\
&=2(x_1x_2+x_3x_4)-2(-x_1x_2-x_3x_4)\\
&=4(x_1x_2+x_3x_4)>0.
\end{aligned} $$ Hence, $\lambda_1(\Gamma^{\ast})>\lambda_1(\Gamma^{\prime}),$ a contradiction.  Thus, Claim \ref{only1negative} holds.

Without loss of generality, by Claim \ref{only1negative},   we can suppose that $\mathscr{C}=v_1v_2\cdots v_{l-1}v_lv_1$ and $v_1v_2$ is the unique negative edge of $\Gamma^{\prime}.$ Let $d_{\Gamma^{\prime}}(v_i)$ denote the degree of $v_i$ in $\Gamma^{\prime}.$
\begin{claim}\label{fenliangdayu0}
$x_i>0$ for $3\leq i\leq n.$
\end{claim}
Otherwise, assume $x_i=0$ for $3\leq i\leq n.$ By Claims \ref{connected} and \ref{only1negative},  $d_{\Gamma^{\prime}}(v_i)\geq 1$ and all edges incident to $v_i$ are positive.
Based on the following equation,   $$0=\lambda_1(\Gamma^{\prime})x_i=\sum\limits_{v_j\in N_{\Gamma^{\prime}}(v_i)}x_j,$$  we have $x_j=x_i=0$ for any $v_j\in N_{\Gamma^{\prime}}(v_i),$ which contradicts  Claim \ref{1zero}.

Next, without loss of generality, we suppose $x_1\geq x_2\geq0.$ 
\begin{claim}\label{k=3}
$l=3.$
\end{claim}
Otherwise, assume $l\geq5$. We first consider that $l\geq6.$ Then $v_{3}\nsim v_{l-1}$ and $v_{1}\nsim v_{l-1},$  and we can construct a new signed graph $\Gamma^{\ast}$  obtained from $\Gamma^{\prime}$ by adding a positive edge $v_3v_{l-1}.$  Clearly, $\Gamma^{\ast}$ is also unbalanced and $\mathcal{C}_4^-$-free. By Rayleigh principle and Claim \ref{fenliangdayu0}, we obtain that
$$\begin{aligned}
\lambda_1(\Gamma^{\ast})-\lambda_1(\Gamma^{\prime})&\geq x^{\top}A(\Gamma^{\ast})x-x^{\top}A(\Gamma^{\prime})x\\
&=2x_3x_{l-1}>0.\\
\end{aligned}$$
Hence, $\lambda_1(\Gamma^{\ast})>\lambda_1(\Gamma^{\prime}),$ a contradiction. Thus, $l=5,$  $\mathscr{C}=v_1v_2v_3v_4v_5v_1,$ $v_1\nsim v_3,$  $v_2\nsim v_5,$ $v_1\nsim v_4$ and $v_2\nsim v_4.$   Let $W_1=N_{\Gamma^{\prime}}(v_1)\cap N_{\Gamma^{\prime}}(v_5)$ and $W_2=N_{\Gamma^{\prime}}(v_2)\cap N_{\Gamma^{\prime}}(v_3).$ We assert that $W_1\neq \emptyset.$ Otherwise, assume $W_1= \emptyset.$ Then  we can construct a new signed graph $\Gamma^{\ast}$  obtained from $\Gamma^{\prime}$ by adding a positive edge $v_2v_5.$  Clearly, $\Gamma^{\ast}$ is also unbalanced and $\mathcal{C}_4^-$-free. By Rayleigh principle, we obtain that
$$\begin{aligned}
\lambda_1(\Gamma^{\ast})-\lambda_1(\Gamma^{\prime})&\geq x^{\top}A(\Gamma^{\ast})x-x^{\top}A(\Gamma^{\prime})x\\
&=2x_2x_5\geq0.
\end{aligned} $$
If $\lambda_1(\Gamma^{\ast})=\lambda_1(\Gamma^{\prime})$, then $x$ is also an eigenvector of $A(\Gamma^{\ast})$
corresponding to $\lambda_1(\Gamma^{\ast}).$ Based on the following equations,
$$\lambda_1(\Gamma^{\prime})x_{2}=\sum_{v_s\in N_{\Gamma^{\prime}}(v_2)}\sigma^{\prime}(v_sv_2)x_s$$ and
$$  \lambda_1(\Gamma^{\ast})x_{2}=\sum_{v_s\in N_{\Gamma^{\prime}}(v_2)}\sigma^{\prime}(v_sv_2)x_s+x_5,$$
we obtain that $x_5=0,$ which   contradicts  Claim \ref{fenliangdayu0}. Hence, $\lambda_1(\Gamma^{\ast})>\lambda_1(\Gamma^{\prime}),$ a contradiction.  Similarly, $W_2\neq \emptyset.$  Recall that $\Gamma^{\prime}$ is $\mathcal{C}_{4}^-$-free. Then for any $v_p\in W_1,$ $v_p\nsim v_2$ and  $v_p\nsim v_3$ and for any $v_q\in W_2,$ $v_q\nsim v_1$ and  $v_q\nsim v_5.$ If $x_5\geq x_3,$ then for all  $v_q\in W_2,$ we can construct a new signed graph $\Gamma^{\ast}$  obtained from $\Gamma^{\prime}$ by rotating all positive edges $v_2v_q$ to  the  non-edge position $v_1v_q,$ rotating all  positive edges $v_3v_q$ to  the  non-edge position $v_5v_q$ and adding a positive edge $v_1v_3.$  Clearly, $\Gamma^{\ast}$ is also unbalanced and $\mathcal{C}_4^-$-free. By Rayleigh principle  and Claim \ref{fenliangdayu0}, we obtain that
 $$\begin{aligned}
\lambda_1(\Gamma^{\ast})-\lambda_1(\Gamma^{\prime})&\geq x^{\top}A(\Gamma^{\ast})x-x^{\top}A(\Gamma^{\prime})x\\
&=2\sum\limits_{v_q\in W_2}x_q(x_1-x_2)+2\sum\limits_{v_q\in W_2}x_q(x_5-x_3)+2x_1x_3\\
&>0.
\end{aligned} $$
Hence, $\lambda_1(\Gamma^{\ast})>\lambda_1(\Gamma^{\prime}),$ a contradiction.  If $x_5 \leq x_3,$ then  for all $v_p\in W_1$ and $v_q\in W_2,$ we can construct a new signed graph $\Gamma^{\ast}$  obtained from $\Gamma^{\prime}$ by rotating all  positive edges $v_5v_p$ to  the  non-edge position $v_3v_p,$ rotating all positive edges $v_2v_q$ to  the  non-edge position $v_1v_q,$  deleting the positive edge $v_2v_3$ and adding two positive edges $v_1v_3$ and $v_2v_5.$  Clearly, $\Gamma^{\ast}$ is also unbalanced and $\mathcal{C}_4^-$-free. By Rayleigh principle  and Claim \ref{fenliangdayu0}, we obtain that
$$\begin{aligned}
\lambda_1(\Gamma^{\ast})-\lambda_1(\Gamma^{\prime})&\geq x^{\top}A(\Gamma^{\ast})x-x^{\top}A(\Gamma^{\prime})x\\
&=2\sum\limits_{v_p\in W_1}x_p(x_3-x_5)+2\sum\limits_{v_q\in W_2}x_q(x_1-x_2)+2x_3(x_1-x_2)+2x_2x_5 \\
&>0.
\end{aligned} $$
Hence, $\lambda_1(\Gamma^{\ast})>\lambda_1(\Gamma^{\prime}),$  a contradiction. Thus, Claim \ref{k=3} holds.

By Claims \ref{containallnegative}, \ref{only1negative} and \ref{k=3}, we have $\mathscr{C}=v_1v_2v_3v_1$ and $v_1v_2$ is the unique negative edge of $\Gamma^{\prime}$.  Without loss of generality, suppose $x_r=\max\limits_{1\leq i\leq n}{x_i}.$
\begin{claim}\label{duwei N-2}
$d_{\Gamma^{\prime}}(v_r)\geq n-2.$
\end{claim}
Otherwise, assume $d_{\Gamma^{\prime}}(v_r)\leq n-3$. Then
$$\lambda_1(\Gamma^{\prime})x_r=\sum_{v_j\in N_{\Gamma^{\prime}}(v_r)}\sigma^{\prime}(v_rv_j){x_j}\leq(n-3)x_r.$$ Thus, $\lambda_1(\Gamma^{\prime})\leq n-3,$
 a contradiction. 

For $S\subset V(G),$ we denote by  $G[S]$  the subgraph of $G$ induced by $S$ and $\Gamma^{\prime}[S]$  the signed induced  subgraph of $\Gamma^{\prime}=(G,\sigma^{\prime})$ whose underlying graph is $G[S]$ and edges have the same signs as them in $\Gamma^{\prime}$. 
\begin{claim}\label{d=n-1}
 $d_{\Gamma^{\prime}}(v_r)=n-1$.
\end{claim}
Otherwise,  assume $d_{\Gamma^{\prime}}(v_r)=n-2$ by Claim \ref{duwei N-2}. We first consider that $r=1.$ Without loss of generality, let  $N_{\Gamma^{\prime}}(v_r)=\{v_2,v_3,\ldots,v_{n-1}\}.$ Then
$$\begin{aligned}\lambda_1(\Gamma^{\prime})x_{r}&=\lambda_1(\Gamma^{\prime})x_{1}=-x_2+\sum_{3\leq s\leq n-1}x_s\\
&\leq 0+(n-3)x_r=(n-3)x_r.
\end{aligned}$$
Hence, $\lambda_1(\Gamma^{\prime})\leq n-3,$ a contradiction.  Similarly, $r\neq2.$ Next, we consider that  $r=3.$ Without loss of generality, let  $N_{\Gamma^{\prime}}(v_r)=\{v_1,v_2,v_5,\ldots,v_{n}\}.$ Since  $\Gamma^{\prime}$ is $\mathcal{C}_4^-$-free, we have $v_i\nsim v_j$ for any $1\leq i\leq2$ and $5\leq j\leq n.$  We assert that  $\Gamma^{\prime}[V(\Gamma^{\prime})\backslash\{v_1,v_2,v_3,v_4\}]=(K_{n-4},+).$  Otherwise,  without loss of generality, assume $v_5\nsim v_6.$  Then we can construct a new signed graph $\Gamma^{\ast}$  obtained from $\Gamma^{\prime}$ by adding a positive edge $v_5v_6.$   Clearly, $\Gamma^{\ast}$ is also unbalanced and $\mathcal{C}_4^-$-free. By Rayleigh principle  and Claim \ref{fenliangdayu0}, we obtain that
$$\begin{aligned}
\lambda_1(\Gamma^{\ast})-\lambda_1(\Gamma^{\prime})&\geq x^{\top}A(\Gamma^{\ast})x-x^{\top}A(\Gamma^{\prime})x\\
&=2x_5x_6>0.
\end{aligned} $$
Hence, $\lambda_1(\Gamma^{\ast})>\lambda_1(\Gamma^{\prime}),$  a contradiction.  Similarly,   $v_4v_i$    is  a positive edge of $\Gamma^{\prime}$ for any $v_i\in V(\Gamma^{\prime})\backslash\{v_3,v_4\}.$ Thus,  $\Gamma^{\prime}=\Gamma_{2}$ (see Figure \ref{fig2}).  However, by Lemma \ref{dayu n3}, $\lambda_1(\Gamma^{\prime})=\lambda_1(\Gamma_{2})<\lambda_1(\Gamma_{1}),$  a contradiction.  Now, there is only one case, i.e.,  $r>3.$   We assert that $v_r\nsim v_3.$ Otherwise,  assume  $v_r\sim v_3.$   Since $d_{\Gamma^{\prime}}(v_r)=n-2,$   at least one of $v_1$ and $v_2$ is adjacent to  $v_r.$  Then there must exist a negative cycle $v_1v_2v_3v_rv_1$ or $v_1v_2v_rv_3v_1$ of order $4,$  a contradiction. Thus,   $N_{\Gamma^{\prime}}(v_r)=V(\Gamma)\backslash\{v_3,v_r\}.$  By similar arguments,  
we have $\Gamma^{\prime}=\Gamma_{2},$  a contradiction. So Claim \ref{d=n-1} holds. 
\begin{claim}\label{r=3}
 $r=3.$
\end{claim}
Otherwise, assume $r\neq 3.$ Recall that $x_r=\max\limits_{1\leq i\leq n}{x_i}$ and $d_{\Gamma^{\prime}}(v_r)=n-1.$  If $r>3,$ then there must exist a negative cycle $v_1v_2v_3v_rv_1$ of order $4$, a contradiction.  If $r=1,$ then we assert that  $x_2=\min\limits_{1\leq i\leq n}{x_i}.$ Otherwise,  assume $x_j=\min\limits_{1\leq i\leq n}{x_i}$ and $j\neq 2.$  Then 
$$\begin{aligned}\lambda_1(\Gamma^{\prime})x_{r}&=\lambda_1(\Gamma^{\prime})x_{1}=-x_2+\sum_{3\leq s\leq n}x_s\\
&\leq (x_j-x_2)+(n-3)x_r\\
&< 0+(n-3)x_r=(n-3)x_r.
\end{aligned}$$
Hence, $\lambda_1(\Gamma^{\prime})< n-3,$ a contradiction. Now, we assert that at least one of $V(\Gamma)\backslash\{v_1,v_2,v_3\}$ is not adjacent to $v_2.$ Otherwise, assume that  all vertices of $V(\Gamma)\backslash\{v_1,v_2,v_3\}$ are  adjacent to   $v_2$.  Based on the following equations,
$$\lambda_1(\Gamma^{\prime})x_1=-x_2+\sum\limits_{i=3}^{n}x_i$$ and
$$\lambda_1(\Gamma^{\prime})x_2=-x_1+\sum\limits_{i=3}^{n}x_i,$$  we get that  $x_1=x_2,$ i.e.,  $\max\limits_{1\leq i\leq n}{x_i}=\min\limits_{1\leq i\leq n}{x_i}.$ Recall that  $\Gamma^{\prime}$ is $\mathcal{C}_4^-$-free. Then for any $4\leq k\leq n,$   $v_k\nsim v_3.$ Thus,  
    $$\lambda_1(\Gamma^{\prime})x_3=x_1+x_2,\ \text{i.e.,} \ \lambda_1(\Gamma^{\prime})x_1=2x_1,$$
 a contradiction. Therefore,  without loss of generality, we can suppose  that  $v_4$ is not adjacent to $v_2.$  
 For any $5\leq i\leq n,$  if $v_i\sim v_2,$  then  $v_i\nsim v_4$ since   $\Gamma^{\prime}$ is $\mathcal{C}_4^-$-free.  We can construct a new signed graph $\Gamma^{\ast}$  obtained from $\Gamma^{\prime}$ by rotating the positive edge $v_iv_2$  to the non-edge position $v_iv_4.$   Clearly, $\Gamma^{\ast}$ is also unbalanced and $\mathcal{C}_4^-$-free. By Rayleigh principle, we obtain that
$$\begin{aligned}
\lambda_1(\Gamma^{\ast})-\lambda_1(\Gamma^{\prime})&\geq x^{\top}A(\Gamma^{\ast})x-x^{\top}A(\Gamma^{\prime})x\\
&=2x_i(x_4-x_2)\geq0.
\end{aligned} $$
If $\lambda_1(\Gamma^{\ast})=\lambda_1(\Gamma^{\prime})$, then  $x$ is also an eigenvector of $\Gamma^{\ast}$
corresponding to $\lambda_1(\Gamma^{\ast}).$ Based on the following equations,
$$\lambda_1(\Gamma^{\prime})x_{4}=\sum_{v_s\in N_{\Gamma^{\prime}}(v_4)}\sigma^{\prime}(v_sv_t)x_s$$ and
$$  \lambda_1(\Gamma^{\ast})x_{4}=\sum_{v_s\in N_{\Gamma^{\prime}}(v_4)}\sigma^{\prime}(v_sv_t)x_s+x_i,$$
we obtain that $x_i=0,$ which  contradicts  Claim \ref{fenliangdayu0}. Hence, $\lambda_1(\Gamma^{\ast})>\lambda_1(\Gamma^{\prime}),$ a contradiction.
Thus,  for any $4\leq k\leq n,$ $v_k\nsim v_2$ and $v_k\nsim v_3.$  Based on the following equations,  $$\lambda_1(\Gamma^{\prime})x_2=x_3-x_1$$ and $$\lambda_1(\Gamma^{\prime})x_3=x_1+x_2,$$ we obtain that
$$\lambda_1(\Gamma^{\prime})x_3=x_3,$$ a contradiction.   Therefore,  $r\neq 1.$ Similarly, $r\neq2.$  So  Claim \ref{r=3} holds.

 By Claims \ref{d=n-1} and \ref{r=3}, we have $d_{\Gamma^{\prime}}(v_3)=n-1.$  
 Thus, $v_i\nsim v_j$ for any $1\leq i\leq2$ and $4\leq j\leq n.$  

\begin{claim}
$\Gamma^{\prime}[V(\Gamma^{\prime})\backslash\{v_1,v_2\}]=(K_{n-2},+).$
\end{claim}
Otherwise,  without loss of generality, assume $v_4\nsim v_5.$  Then we can construct a new signed graph $\Gamma^{\ast}$  obtained from $\Gamma^{\prime}$ by adding a positive edge $v_4v_5.$   Clearly, $\Gamma^{\ast}$ is also unbalanced and $\mathcal{C}_4^-$-free. By Rayleigh principle  and Claim \ref{fenliangdayu0}, we obtain that
$$\begin{aligned}
\lambda_1(\Gamma^{\ast})-\lambda_1(\Gamma^{\prime})&\geq x^{\top}A(\Gamma^{\ast})x-x^{\top}A(\Gamma^{\prime})x\\
&=2x_4x_5>0.
\end{aligned} $$
Hence, $\lambda_1(\Gamma^{\ast})>\lambda_1(\Gamma^{\prime}),$ a contradiction.

 Above all,   $\Gamma^{\prime}=\Gamma_{1},$  which means $\Gamma$ is switching equivalent to $\Gamma_{1}$. This completes the proof.
\end{proof}

\end{document}